\documentclass[12pt]{article}
\usepackage[dvips]{graphics}
\begin{document}

\title{The Modelling of Degenerate Neck Pinch Singularities in Ricci Flow by Bryant Solitons}

\author{David Garfinkle
\thanks {Email: garfinkl@oakland.edu}\\
Department of Physics, Oakland University\\
Rochester, Michigan 48309\\
James Isenberg
\thanks {Email: jim@newton.uoregon.edu} \\
Department of Mathematics, University of Oregon\\
Eugene, OR}

\maketitle

\begin{abstract}
In earlier work, carrying out numerical simulations of the Ricci flows of 
families of rotationally symmetric geometries on $S^3$, we have found strong 
support for the contention that (at least in the rotationally symmetric case) 
the Ricci flow for a ``critical" initial geometry-- one which is at the 
transition point between initial geometries (on $S^3$) whose volume-normalized 
Ricci flows develop a singular neck pinch, and other initial geometries whose 
volume-normalized Ricci flows converge to a round sphere--evolves into a 
``degenerate neck pinch". That is, we have seen in this earlier work that 
the Ricci flows for the critical geometries become locally cylindrical in 
a neighborhood of the initial pinching, and have the maximum amount of
curvature at one or both of the poles. Here, we explore the behavior of 
these flows at the poles, and find  strong support for the conjecture that 
the Bryant steady solitons accurately model this polar flow.

\end{abstract}
\section{Introduction}

While a considerable amount has been learned during the past five 
years \cite{hamilton3, Per1, Per2} \cite{CLN, RFV2} about the Ricci flow of three 
dimensional Riemannian geometries,  and while this knowledge has been used 
very successfully to study the relationship between topology and 
geometry \cite{Per1, CZ, MT, KL}, there are still many unanswered questions 
regarding the details of three dimensional Ricci flow, especially concerning 
flows which develop singularities. One of these questions concerns the
circumstances and the details of the formation of singularities of Type II. 
The Type II singularities are characterized by the condition that, if the 
singularity occurs at a finite time $T$, then the quantity $|Rm(t)|(T-t)|$ 
blows up as $t$ approaches $T$; this behavior contrasts with that of Type 
I singularities, for which this same quantity is uniformly bounded. 

It is believed that the standard neck pinch singularities, which according 
to Hamilton's scenario  \cite{HamForm} play a major role in three dimensional 
Ricci flow dynamics, are Type I. The work of Angenent and 
Knopf \cite{AK1, AK2} supports this contention, and describes some of the 
detailed asymptotic behavior of standard neck pinch singularities, at least 
in the rotationally symmetric case. 

It has been conjectured \cite{HamForm}, it has been demonstrated 
numerically \cite{numericalricci}, and it has now been 
proven \cite{GuZhu}, that Type 
II singularities occur during the course of a Ricci flow which uses 
a ``critical geometry'' for initial data. To obtain a critical geometry in 
the sense we mean here, one considers a one-parameter family of 
(initial data) Riemannian metrics on $S^3$, with the volume-normalized
Ricci flow developing standard neck pinch singularities for low values of 
the parameter, and with singularities avoided for those flows which start 
from geometries with high values of the parameter. The critical geometry 
then corresponds to the transitional (boundary) value of the parameter. 

Our earlier  work \cite{numericalricci}, in which we have 
numerically simulated the 
Ricci flow of critical geometries, uses families of geometries which are 
all rotationally symmetric and also reflection symmetric across the equator. 
The geometries are all ``corsetted spheres'', with the parameter measuring 
the degree of corsetting relative to a round sphere. For tight corsetting 
(small parameter) the Ricci flow develops a standard neck pinch singularity, 
as described in \cite{AK1}. For loose corsetting, the flow approaches the 
round sphere. In \cite{numericalricci}, we numerically 
simulate the Ricci flow of both 
tightly and loosely corsetted geometries, and we are then able to focus 
on the transitional critical geometry, and examine its flow. 

The asymptotic behavior of the Ricci flow of the critical geometry seen in our 
simulations is quite different from that which we see starting at geometries 
with non critical values of the parameter. The curvature concentrates at 
the poles (recall our assumption of reflection symmetry), and overall the 
geometry approaches that of a ``javelin'': it becomes increasingly 
cylindrical everywhere except at the poles, where the curvature blows up. 
In analogy with similar behavior seen for critical mean curvature flows, 
this has been labeled a ``degenerate neck pinch". 

In our earlier work, we did not closely explore the details of the Ricci 
flow for our critical geometries at the poles. We noted the curvature blow 
up, but nothing further. Here, we focus on the flow at the poles, and we 
verify that, as has been conjectured, the flow in the neighborhood of the 
poles is very accurately modelled locally by the flow of the Bryant steady
soliton \cite{Bryant}; see also chapter 1 of \cite{RFV2}.

The Bryant steady soliton is the unique (up to scaling) rotationally symmetric Ricci 
gradient soliton on $R^3$ which neither shrinks nor expands. That is, the 
metric $g$ on $R^3$ is 
rotationally symmetric about a fixed point on $R^3$, taking the form 
\begin{equation}
g=dr^2 + a(r)^2 \gamma_{round}
\label{rotsymmetric}
\end{equation}
for the round metric $\gamma_{round}$ on $S^2$ and for a positive function 
$a(r)$; and it satisfies the steady Ricci gradient soliton equation
\begin{equation}
{R_{ab}}+{\nabla _a}{\nabla _b} f = 0 
\label{Ricsol}
\end{equation} 
for a function $f$. 

If one substitutes into equation (\ref{Ricsol}) the metric 
$g$ of the form (\ref{rotsymmetric}) together with a rotationally symmetric 
function $f(r)$, one obtains a coupled system of ordinary differential
equations, to be solved for $a(r)$ and $f(r)$.  In \cite{Bryant}, 
(See also chapter 1 of \cite{RFV2}) this ODE system is written 
out and analyzed, it is shown that there is a unique solution (up to 
homothety), and the profile of this solution (obtained by numerical 
integration) is exhibited. We use this profile in our model matching here. 
(See section 4.)

One of the key features of a steady Ricci soliton $(g,f)$ is that the Ricci 
flow of the soliton metric $g$  fixes $g$ up to a time-dependent 
diffeomorphism (generated by the vector field $g^{-1} ( df,\_))$. Hence the 
curvature is time independent. The sense in which such a flow can model 
the singularity developing at the pole of our degenerate neck pinch is in 
terms of the blow up (rescaling) of that singularity. Specifically, given a 
Ricci flow
solution $g(t)$ which is going singular at the point $x^*$ at time $T$, 
let us define the 
``blow up'' metrics $\tilde{g}(t)$ by setting
\begin{equation}
\label{blowup}
\tilde{g} ( t) :=  \rho(t) g( t )
\end{equation}
where the function $\rho (t)$ 
is chosen so that all the rescaled metrics $\tilde{g} (t)$ have the 
same value of $|Rm|$ at the point $x^*$.
The 
singularity is then modeled by the soliton if, for spatial points near
$x^*$ and
for $t$ approaching $T$, the blow up metrics $\tilde{g}$ approach closer 
and closer to the soliton metric. We provide strong numerical evidence
here that 
for the degenerate neck pinch solutions which we study here, the 
geometry at the 
poles is indeed modelled by Bryant steady solitons in this sense. 

Since, in any family of initial geometries,  it is difficult to precisely find the 
critical initial geometry which flows (via Ricci flow) 
to a degenerate neck pinch, for our numerical studies we have chosen to carry out 
a slightly different comparison as well: we consider a sequence of initial geometries $g_{\alpha}$
which are all sub critical (loose corsetting) and which approach the critical geometry. Since these 
are sub critical, each of  their flows evolve toward a time $t_{\alpha}$ (different for each one)
at which the curvature at the poles reaches a maximum, after which their flows dissipate
the curvature. So one alternative way to test the Bryant steady soliton modelling conjecture 
is to evolve each of the $g_{\alpha}$ geometries up to its time of maximum curvature $t_{\alpha}$
and then compare the geometry near the pole of $g_{\alpha} (t_{\alpha} )$ with 
the Bryant steady soliton, scaled to have the same maximum geometry as  $g_{\alpha} (t_{\alpha} )$. Our numerical work indicates that this comparison too matches very well. 

When our initial work \cite{numericalricci}  on the numerical simulation 
of the Ricci 
flow of degenerate neck pinches was carried out, it had not been shown
mathematically that Type II singularities do develop during Ricci flow. 
This has now been shown by Gu and Zhu \cite{GuZhu}. Their work does not, 
however, tell us any of the details of such  Ricci flows. Our work here 
indicates what those details should be. One hopes that a proof of this
behavior will be forthcoming.  

\section{Ricci flow}

The types of metrics considered are the same as in \cite{numericalricci}. 
We have spherically symmetric metrics on $S^3$ which take the form
\begin{equation}
d {s^2} = {e^{2X}}  \left ( {e^{- 2 W}}  d {\psi ^2}  +  {e^{2 W}}
{\sin ^2} \psi  [ d {\theta ^2}  +  {\sin ^2} \theta  d {\phi ^2} ].
\right )
\label{metric}
\end{equation}
The corseted sphere geometries that we use for initial data have 
$W=X$ with $X$ determined intrinsically by
\begin{eqnarray}
4{e^{4X}} {\sin^2} \psi &=& {\sin^2} 2 \psi \; \; \; {\rm for} \; 
{\cos^2}\psi \ge {\textstyle {\frac 1 2}} 
\nonumber
\\
4{e^{4X}} {\sin^2} \psi &=& {\sin^2} 2 \psi + 4 \lambda {\cos^2} 2 \psi 
\; \; \;  {\rm for} \;
{\cos^2}\psi \le {\textstyle {\frac 1 2}},
\end{eqnarray}
where $\lambda $ is a constant that determines the amount of corsetting.
In particular $\lambda =0$ corresponds to two round three spheres joined 
at the poles, while for positive $\lambda$ this cusp is smoothed out.

To ensure that  the evolution is well behaved numerically, we evolve using the
volume-normalized DeTurck flow \cite{DeTurck}
\begin{equation}
{\partial_t}{g_{ab}} = - 2 {R_{ab}} + 2 {D_{(a}}{V_{b)}} + 
{\textstyle {\frac 2 3}} {\hat r} {g_{ab}}
\label{DeT}
\end{equation}
Here the spatial constant $\hat r$ is the volume average of the scalar curvature,
and the vector field $V^a$ is given by
\begin{equation}
{V^a} = {g^{bc}}\left ( {\Gamma ^a _{bc}} - {\Delta ^a _{bc}} \right ),
\end{equation}
where $\Gamma ^a _{bc}$ is the connection of the metric $g_{ab}$ while
$\Delta ^a _{bc}$ is the connection of a round three sphere.  This evolution equation
for the metric $g$, together with the form of $g$
in (\ref{metric}), yields partial 
differential equations for the metric quantities $X$ and $W$.  However,
we find that it is more convenient to use $X$ and the quantity 
$S \equiv W/{\sin^2} \psi$.  The reason for this is that smoothness of
the metric requires that $W$ vanish at the poles of the three sphere at the rate of
${\sin ^2}\psi$ and this condition is automatically enforced by smoothness
of $S$.  The evolution equations  for $X$ and $S$ that follow from 
equation (\ref{DeT}) are 
\begin{eqnarray}
{\partial _t} X = {e^{2(W-X)}}  \bigg [ {X''}  +  2  \cot \psi
{X' }  -  2  +  {\frac 1 2}  ( {{[{X'}]} ^2}  +  {{[{W'}]} ^2} )
+  3
{X'}  {W'}
\nonumber
\\
+  (1 \, - \, {e^{- 4 W}})  \left ( 
{\frac 1 {2 {\sin ^2} \psi }}  +  1  +  2  \cot \psi  {W'} \right )  \bigg ]  + 
{\frac {\hat r}  3}.
\label{evolveX}
\\
{\partial _t} S = {e^{2(W-X)}} \biggl [ {S''} + 6 \cot \psi {S'}
- 8 S  - {\frac 3 {2 {\sin ^4} \psi }}
\left ( 1 - 4 W - {e^{-4W}}\right )
\nonumber \\
+ {\frac {1 - {e^{-4W}}} {{\sin ^2} \psi }}
\left ( 1 - 2 [\cot \psi {X'} + 2 \sin \psi \cos \psi {S'} + 4
{\cos ^2} \psi S]\right )
\nonumber \\
- {\frac 1 2} \bigl ( {{[{X'}/\sin \psi ]}^2} +
{{[\sin \psi {S'}+2 \cos \psi S]}^2} 
\nonumber \\
+ 6 [{X'}/\sin \psi ] [ \sin \psi
{S'} + 2 \cos \psi S] \bigr ) \biggr ].
\label{evolveS}
\end{eqnarray}  
Here a prime denotes differentiation with respect to $\psi$ and the quantities $W$ 
and $W'$ should be thought of as derived from $S$ through 
$W = S {\sin ^2} \psi$ and 
${W'} = {S'} {\sin ^2} \psi + 2 S \sin \psi \cos \psi$.   
The quantity $\hat r$ is calculated by
\begin{equation}
{\hat r} = {2 \over N} \; {\int _0 ^\pi }  d \psi  {e^{X + 3 W}}  
\left ({e^{ - 4 W}}  - 1  -  4  \sin \psi  \cos \psi {W'}  +  
{\sin^2}
\psi [ 3 + {{({X'} + {W'})}^2} ] \right ),
\end{equation}
where the normalization constant $N$ is given by
\begin{equation}
N \equiv {\int _0 ^\pi } \; d \psi \; {e^{3X + W}} \; {\sin ^2} \psi.
\end{equation}

\section{Numerical methods}
To numerically integrate the PDEs for $X$ and $S$, we approximate these 
functions by their values on a grid and approximate the PDEs by 
finite difference equations.  Let $F$ stand for the pair $(X,S)$
and define the numbers $F^n _i$ by 
\begin{equation}
{F^n _i} = F((i-{\textstyle {\frac 3 2}})\Delta \psi ,n\Delta t) 
\end{equation}   
That is, the $F^n _i$ are the values of $F$ on a grid with spatial
step $\Delta \psi$ and time step $\Delta t$.  For a PDE of the form
${\partial _t} F = {\cal O} F $ for some operator $\cal O$ we
approximate $\cal O$ by the finite difference operator ${\hat {\cal O}}$
given by replacing all spatial derivatives in $\cal O$ with centered
finite differences.  
We then approximate ${\partial _t}F$ by $({F^{n+1} _i} - {F^n _i})/\Delta t$.
Since at any given time step we know $F^n _i$ and want to solve for
$F^{n+1} _i$ the simplest thing to do is to apply $\hat {\cal O}$ at
time step $n$, which yields the finite difference equation
\begin{equation}
{\frac {{F^{n+1} _i} - {F^n _i}} {\Delta t}} = {\hat {\cal O}}{F^n _i},
\end{equation}
which has the solution 
\begin{equation}
{F^{n+1} _i} = {F^n _i} + \Delta t {\hat {\cal O}}{F^n _i}.
\end{equation} 
This is the method used in the work described in  \cite{numericalricci}.  Unfortunately,
this method is quite slow for the following reason: as is typical for
parabolic equations, stability of the numerical method requires a time
step $\Delta t$ that is of order $({\Delta \psi})^2$. This means that
it becomes extremely slow to run simulations with high resolution.  For
this reason--limitations in  resolution--we were able in \cite{numericalricci} 
to show the existence of a critical solution, but we were not able to examine
that critical solution accurately enough to characterize it.  To obtain the 
needed level of accuracy without sacrificing efficiency,
we need a faster numerical method--one that does not have such a limitted
time step.  One way to overcome the limitation on the time step is to 
apply the operator $\hat {\cal O}$ at time step $n+1$ rather than time step
$n$.  This yields 
\begin{equation}
{F^{n+1} _i} = {F^n _i} + \Delta t {\hat {\cal O}}{F^{n+1} _i}.
\end{equation}
At first this equation does not seem helpful, since we are given $F^n _i$ 
and we want to find $F^{n+1} _i$.  However, this equation has the solution
\begin{equation}
{F^{n+1} _i} = {{[I-\Delta t {\hat {\cal O}}]}^{-1}} ({F^n _i}),
\end{equation}
where $I$ is the identity operator. This equation has the advantage that 
stability places no restriction on the size of the time step.  It also has the 
disadvantage that it requires the inversion of the nonlinear
operator $I-\Delta t {\hat {\cal O}}$ a process that is both difficult
and slow.  However, it is only the principal part of the operator 
$\hat {\cal O}$ that leads to the restriction to small time step.  The
solution to this dilemma is then to split $\hat {\cal O}$ into two parts,
one of which contains the principal part but is also simple enough to 
be inverted quickly.  We then apply that part at time step $n+1$ and
the rest at time step $n$.  Specifically, we write the equations of
motion as 
\begin{eqnarray}
{{\hat {\cal O}}_1}({X^{n+1} _i}) 
=  {{\hat {\cal O}}_3} ({X^n _i},{S^n _i}) 
\\
{{\hat {\cal O}}_2}({S^{n+1} _i})
=  {{\hat {\cal O}}_4} ({X^n _i},{S^n _i}),
\end{eqnarray}
where the operators above are the finite difference version of 
\begin{eqnarray}
{{\cal O}_1} (X) = {\frac {e^{2(X-W)}} {\Delta t}} X - ({X''} + 2 
\cot \psi {X'})
\\
{{\cal O}_2} (S) = {\frac {e^{2(X-W)}} {\Delta t}} S - ({S''} + 6 
\cot \psi {S'}) 
\\
{{\cal O}_3} (X,S) = {\frac {e^{2(X-W)}} {\Delta t}} X
 -  2  +  {\frac 1 2}  ( {{[{X'}]} ^2}  +  {{[{W'}]} ^2} )
+ 3 {X'}  {W'}
\nonumber
\\
+  (1 \, - \, {e^{- 4 W}})  \left (
{\frac 1 {2 {\sin ^2} \psi }}  +  1  +  2  \cot \psi  {W'} \right )   
+ {\frac {\hat r} 3} {e^{2(X-W)}}
\\
{{\cal O}_4} (X,S) = {\frac {e^{2(X-W)}} {\Delta t}} S 
- 8 S  - {\frac 3 {2 {\sin ^4} \psi }}
\left ( 1 - 4 W - {e^{-4W}}\right )
\nonumber \\
+ {\frac {1 - {e^{-4W}}} {{\sin ^2} \psi }}
\left ( 1 - 2 [\cot \psi {X'} + 2 \sin \psi \cos \psi {S'} + 4
{\cos ^2} \psi S]\right )
\nonumber \\
- {\frac 1 2} \bigl ( {{[{X'}/\sin \psi ]}^2} +
{{[\sin \psi {S'}+2 \cos \psi S]}^2}
\nonumber \\
+ 6 [{X'}/\sin \psi ] [ \sin \psi
{S'} + 2 \cos \psi S] \bigr ) .
\end{eqnarray}
In the expressions for ${\cal O}_1$ and ${\cal O}_2$ the quantity 
$e^{2(X-W)}$ is evaluated at time step $n$ even though the operator 
is applied to an argument at time step $n+1$.  
Thus given the metric functions $(X,S)$ at time step $n$ we produce
their values at time step $n+1$ through
\begin{eqnarray}
{X^{n+1} _i}= {{{\hat {\cal O}}_1}^{-1}}
({{\hat {\cal O}}_3} ({X^n _i},{S^n _i}))
\\
{S^{n+1} _i}= {{{\hat {\cal O}}_2}^{-1}}
({{\hat {\cal O}}_4} ({X^n _i},{S^n _i})).
\end{eqnarray}
The operators ${{\hat {\cal O}}_1}$ and ${{\hat {\cal O}}_2}$ are 
linear operators that are easily and rapidly inverted using 
the cyclic tridiagonal method as given in \cite{sauletal}

\section{Results and comparison with the Bryant steady soliton}
 
\begin{figure}
\includegraphics{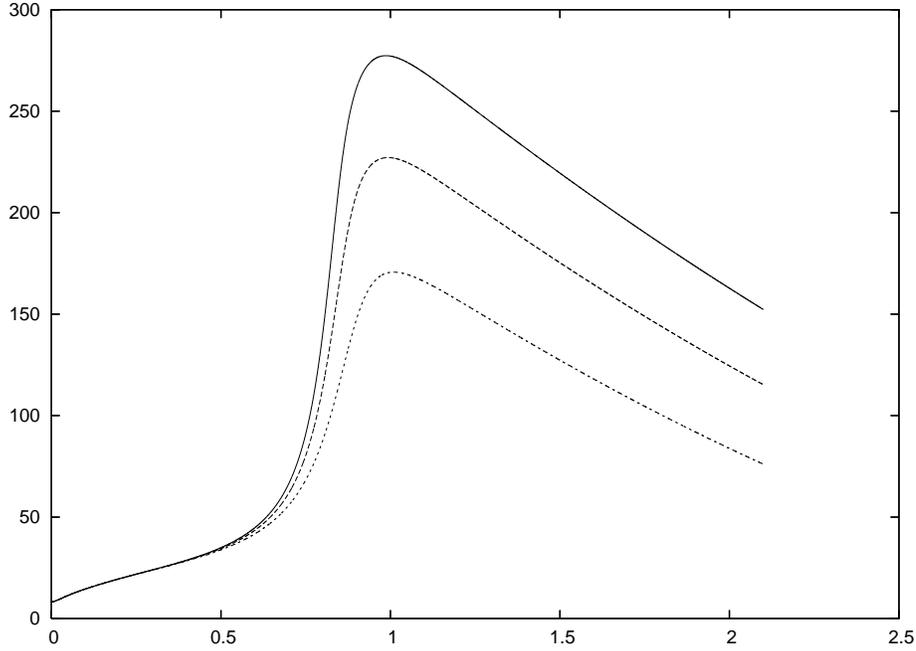}
\caption{\label{fig1}
Plot of scalar curvature at the pole {\it vs} time 
for three different values of $\lambda$}
\end{figure}

All runs of the computer code have been  carried out in double precision 
with 10,000 spatial grid points and with $dt=d\psi$.  
Through a binary search, we have determined
the critical value of $\lambda$ and have then examined the behavior of
the curvature for several slightly subcritical solutions.  In all
cases, during the course of Ricci flow, the curvature becomes  large and 
then diminishes, with the
maximum of the curvature occuring at the poles.  Figure (\ref{fig1}) shows 
the scalar curvature at the pole as a function of time for three
different subcritical solutions.  

It is clear from this figure that 
as the initial data gets closer to the critical data, the maximum 
curvature occurring during the corresponding flow gets larger.  This suggests
that  to get the best view (via numerical simulation) of the critical 
solution, one should examine the evolving  metric at the time for which the 
curvature of a slightly subcritical solution is at its maximum.  

As noted in the introduction, we follow two approaches here to see if 
the behavior of the critical solution at the poles is modeled accurately by 
the Bryant steady soliton. Before carrying either of them through, we need a 
numerical simulation of the metric for this soliton. To obtain this, we first substitute
a metric of the form (\ref{rotsymmetric}) into the Ricci soliton 
equation (\ref{Ricsol}), and we are led to the following ODE system

\begin{eqnarray}
{f''} &=& {\frac {2 {a''}} a} 
\label{BrySol1}
\\
{a''} &=& {f'} {a'} + {\frac {1 - {{({a'})}^2}} a} ,
\label{BrySol2}
\end{eqnarray} 
where a prime denotes differentiation with respect to $r$. Since this is a 
system of second order ODEs for two functions, one might expect a four
parameter family of solutions.  However, equations (\ref{BrySol1})-(\ref{BrySol2})
are invariant under the addition of a constant to $f$
and smoothness of the solution at the origin imposes
additional constraints.  In the end, the general smooth solution is determined
by a single parameter $\alpha$ which appears in the series expansions for $a$ 
and $f'$ around the origin as follows
\begin{eqnarray}
a &=& r + \alpha {r^3} + \dots
\label{powera}
\\
{f'} &=& 12 \alpha r + \dots
\label{powerf}
\end{eqnarray}
where the dots $\dots $ stand for higher order terms.  However, even this 
degree of freedom is to a certain extent misleading.  It is only for
negative $\alpha$ that nonsingular solutions exist for $0\le r < \infty$
and solutions with different negative values of $\alpha$ differ only
in an overall constant scale factor in the metric.   
Thus up to scale, there is really only one Bryant steady soliton.  The relation between
the overall scale and the parameter $\alpha$ is reflected in the expression
 $R=-36 \alpha$ for the scalar curvature at the origin.
 
 To develop the Bryant steady soliton using (\ref{BrySol1}) and (\ref{BrySol2}), we first 
 choose $\alpha$ (the choice depending on the desired scale; see below) and use it to 
 express the expansion forms for $a$ and $f'$ in a neighborhood of the origin, following
(\ref{powera}) and (\ref{powerf}).  We then numerically integrate the rest 
of the way using the fourth order Runge-Kutta method.\cite{sauletal}

For a spherically symmetric metric, the Ricci tensor has two independent
eigenvalues: $R_{S^2}$ the eigenvalue in the directions tangent to the
symmetry $S^2$ and $R_\perp$ the eigenvalue in the direction 
perpendicular to the symmetry $S^2$.  For both our Ricci flow simulations
and for our numerical integration of the Bryant soliton ODEs, we 
calculate these Ricci eigenvalues and compare them. 

\begin{figure}
\includegraphics{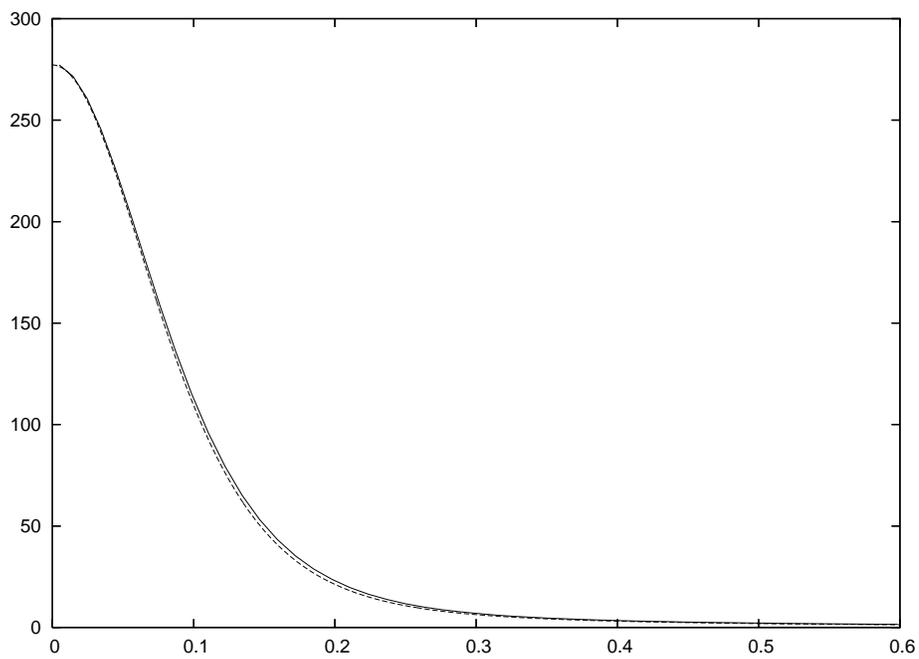}
\caption{\label{r1compare}$R_\perp$ as a function of radial length for a near
critical solution at the time of maximum curvature (solid line) and for the
Bryant soliton (dashed line)}
\end{figure}

\begin{figure}
\includegraphics{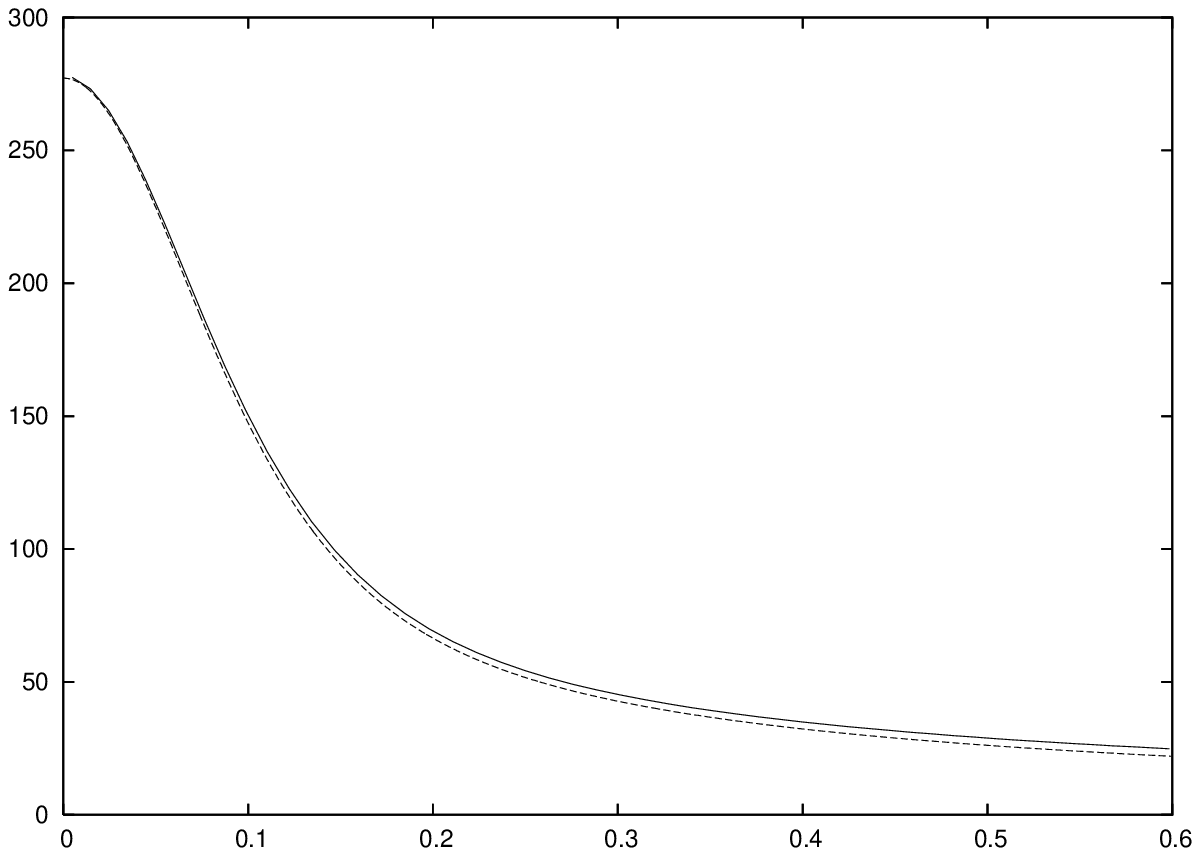}
\caption{\label{r2compare}$R_{S^2}$ as a function of radial length for a near
critical solution at the time of maximum curvature (solid line) and for the
Bryant soliton (dashed line)}
\end{figure}

For one of our comparison studies, we compare the geometry of appropriately scaled 
Bryant steady soliton solutions with  the geometries near the poles at the time of maximum curvature for a sequence of subcritical initial geometries approaching the critical solution.
Here, we choose the parameter $\alpha$ so that the curvature at the tip of the 
soliton matches that at the poles of the flows for the sub critical geometries (at maximum curvature). Figure (\ref{r1compare}) 
gives the $R_\perp$
eigenvalue for the Ricci flow simulation at the time of maximum
curvature (solid line) and for the Bryant soliton (dashed line).  These
quantities are plotted as functions of radial distance.  
Figure (\ref{r2compare}) makes
the same comparison for the $R_{S^2}$ eigenvalue.  Note that in 
both cases there is an excellent match between the Ricci flow simulation
and the Bryant soliton. 

\begin{figure}
\includegraphics{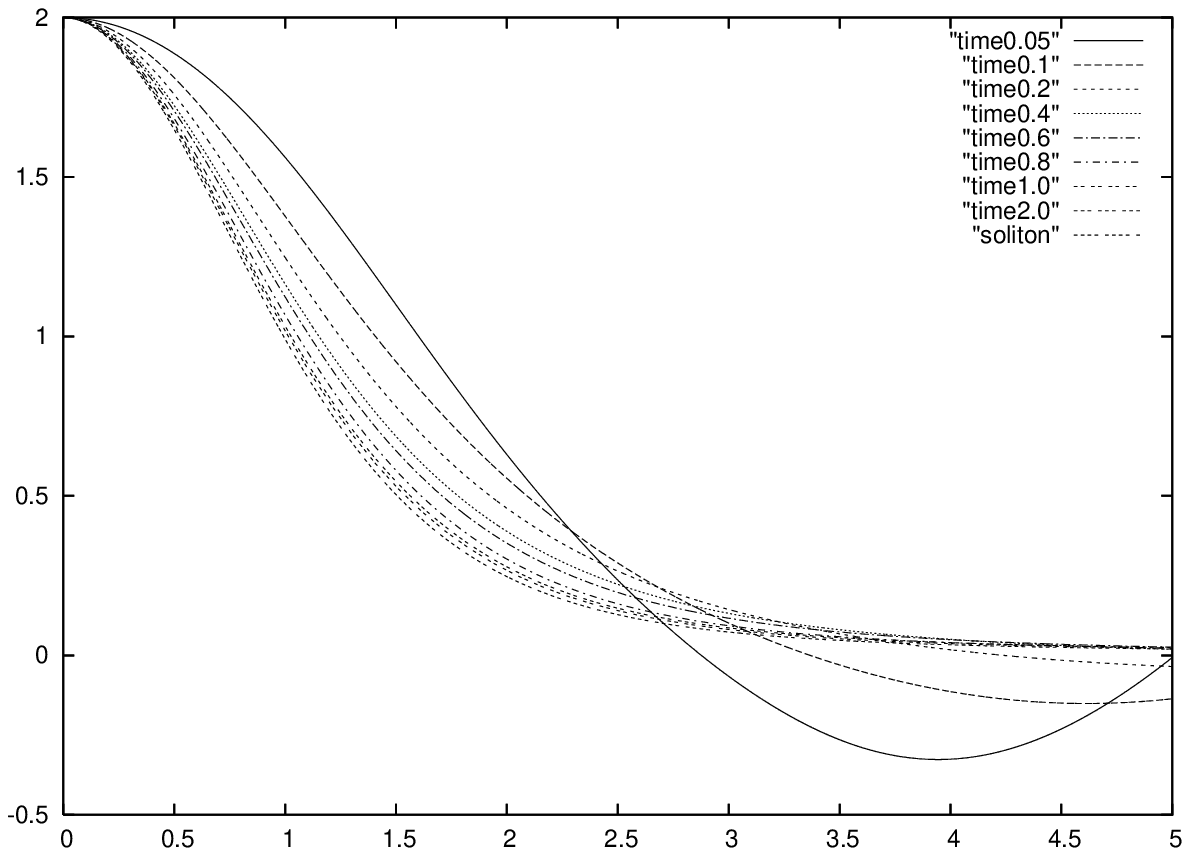}
\caption{\label{rperpscale}
Plot of $R_{\rm perp}$ as a function of radial length both for the 
rescaled metrics of several different times and for the Bryant steady 
soliton} 
\end{figure}

\begin{figure}
\includegraphics{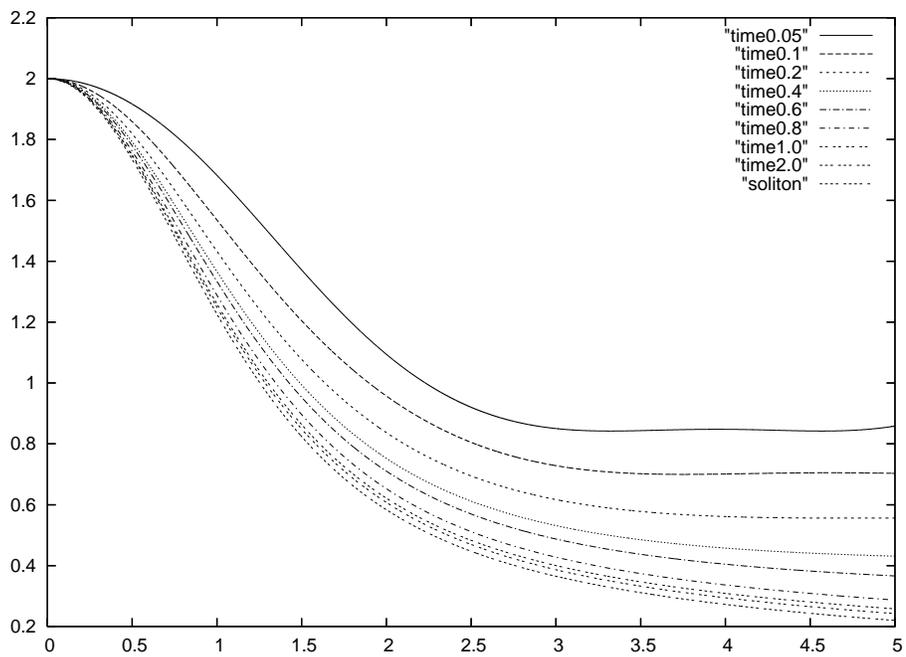}
\caption{\label{rs2scale}
Plot of $R_{S^2}$ as a function of radial length both for the
rescaled metrics of several different times and for the Bryant steady
soliton}
\end{figure}

For our other study, we choose a subcritical initial metric which is 
very close to the 
critical geometry, and consider a sequence of times approaching the 
time of maximum curvature
at the pole. At each time, we calculate the blow up geometry as specified 
in equation 
\ref{blowup}. We scale the blowups, and scale the soliton, so that they 
all have identical curvature at the poles. The results are graphed in 
Figure (\ref{rperpscale}) for $R_{\rm \perp}$ and in Figure (\ref{rs2scale})
for $R_{S^2}$. 
In these simulations, the maximum curvature occurs at a time of approximately
1.0.  We see that the blowup geometries approach agreement with the 
Bryant steady soliton at somewhat early times and retain this agreement for times 
beyond the time of maximum curvature.  
This agreement
(in a neighborhood of the pole) becomes remarkably close. Thus both of our 
tests strongly support the conjecture that the Bryant steady soliton 
accurately models the pole behavior of the Ricci flow of critical geometries 
which develop degenerate singularities.

\section{Conclusions}

Since for each one parameter family of initial geometries one expects degenerate pinches to occur for just a single value of the parameter, direct numerical testing of the behavior of degenerate neck pinches is essentially impossible. Our numerical studies here do, however, strongly support the contention that degenerate neck pinches are modeled very accurately by the Bryant steady soliton, at least in the case of rotationally symmetric geometries.

The natural next step for thse studies is to consider one paramter families of geometries which are \textit{not} rotationally symmetric. Numerical simulation of the Ricci flow for such metrics is expected to be considerably more challenging, and likely will require working on multiple overlapping patches on $S^3$. Carrying out these simulations, however, should allow us to explore whether the Ricci flows of non rotationally symmetric geometries tend to evolve toward rotationally symmetric geometries, both in the case of neck pinch singularity formation, and degenerate neck pinch singularity formation. 

Once numerical evidence for a particular behavior in solutions of a PDE system has been obtained, there is strong motivation to mathematically prove that the behavior is present. As noted above, the existence of Type II singularities in the Ricci flows of critical type geometries has now been proven \cite{GuZhu}. However, the features of these singularities, including the formation of javelin geometries and the modeling by Bryant steady solitons, remains mathematically unverified. This should be a promising direction for future research.

\section{Acknowledgments}

This work was supported by NSF grant PHY-0456655 to Oakland University 
and both PHY-0354659 and PHY-0652903 to The University of
Oregon. We thank Matt Choptuik for helpful discussions.  We also thank both
the University of California at San Diego and the Albert Einstein Institute in Golm, Germany for hospitality while some of this work was carried out.

\end{document}